\documentclass[12pt]{article}
\usepackage{amssymb}
\usepackage{amsmath}

\newtheorem{theo}{Theorem}
\newtheorem{prop}{Proposition}
\newtheorem{lemm}{Lemma}

\newtheorem{rema}{Remark}

\newcommand{\cqfd}
{%
\mbox{}%
\nolinebreak%
\hfill%
\rule{2mm}{2mm}%
\medbreak%
\par%
}
\newfont{\gothic}{eufb10}
\date{\empty}
\begin{document}
\title{On the homotopy types of compact K\"ahler and complex
projective manifolds}
\author{Claire Voisin\\ Institut de math{\'e}matiques de Jussieu, CNRS,UMR
7586} \maketitle \setcounter{section}{-1}
\section{Introduction}
The celebrated Kodaira theorem \cite{kodaira} says that a compact complex manifold is projective
if and only if it admits a K\"ahler form whose cohomology class is integral.
This suggests that K\"ahler geometry is an extension of projective geometry, obtained
by relaxing the integrality condition on a K\"ahler class.
This point of view, together with the
many restrictive conditions on the topology of
K\"ahler manifolds provided by Hodge theory
(the strongest one being the formality theorem \cite{degrimo}),
 would indicate that compact K\"ahler manifolds and complex projective ones
cannot be distinguished by topological invariants.
This is supported by the results known for K\"ahler surfaces, for which
a much stronger statement is known, as a consequence of Kodaira's classification : recall first that two
compact complex manifolds $X$ and $X'$ are said to be deformation equivalent
if there exist a family $\pi:\mathcal{X}\rightarrow B$, where $B$ is a connected analytic space
and $\pi$ is smooth and proper, and two points $b,\,b'\in B$ such that
$$X_b\cong X,\,\,X_{b'}\cong X'.$$
Kodaira shows the following :
\begin{theo} \label{kod}\cite{kodairasurf} A compact K\"ahler surface is deformation equivalent to a
complex  projective
  surface.
\end{theo}
Note that since the basis $B$ of the family is path connected, Ehresmann's theorem shows
that two deformation equivalent complex manifolds are homeomorphic, and
a fortiori  have the same
homotopy type.

Since Kodaira's result, the situation for higher dimensional manifolds remained unknown.
A classical problem, sometimes
called the Kodaira problem (see \cite{dep}), asks whether theorem \ref{kod} still holds in higher dimensions.
Many known examples, like hyperk\"ahler manifolds, or more generally
$K$-trivial K\"ahler manifolds, are deformations of projective ones (see \cite{beauville}).

In this paper, we show however the following result, which in particular provides
a negative answer to Kodaira's problem:
\begin{theo}\label{mainintro} In any dimension
$\geq 4$, there exist compact K\"ahler manifolds which do not have the homotopy type of a
complex projective manifold.
\end{theo}

The case of dimension $3$ remains open.

 Examples are constructed in section \ref{2}. The topological invariant we use here is the
 integral  cohomology
ring, so that a stronger version of the theorem above is the fact that
there exist compact K\"ahler manifolds which do not have a cohomology ring isomorphic to the one
 of a
complex projective manifold. In section \ref{3}, we  explain arguments due to Deligne,
(letter to the author,)
which show how to get similar statements with the rational cohomology ring and, what is still
more surprising, with the complex cohomology algebra.

The examples considered above are birationally equivalent to complex tori, and indeed, the
proofs consist in showing that the cohomology algebra prevents the $H^1$  being endowed
with a polarized Hodge structure.
In section \ref{3},
we construct however a simply connected K\"ahler manifold satisfying the conclusion
of Theorem \ref{mainintro}.  The proof that its rational cohomology ring is not
the one of a projective complex manifold makes a systematic use of Deligne's Lemma
\ref{ledel} (see section \ref{3}), combined with the Hodge index theorem.

\vspace{0.5cm}

{\bf Thanks.} I would like to thank J.-P. Demailly who brought Kodaira's problem to my attention,
 A. Beauville for pointing out an inaccuracy in the first version of this paper, and the
  referees for their careful reading and helpful comments.
I also thank  P. Deligne for his interest in this work and for communicating his insights, which led to
the contents of section \ref{3}.
\section{\label{1} Complex tori with endomorphisms  }
Let $\Gamma \cong\mathbb{Z}^{2n}$ be a lattice, and let $\phi:\Gamma\rightarrow\Gamma$ be an endomorphism.
Let $$\Gamma_\mathbb{C}:=\Gamma\otimes\mathbb{C},$$
on which the $ \mathbb{C}$-linear extension $\phi_\mathbb{C}$ of $\phi$ acts. We shall assume the following
properties :

 (P) {\it The eigenvalues of $\phi$ all have multiplicity $1$, and none of them is real.}

 We choose, among the $2n$ eigenvalues of $\phi$, a set of $n$ elements,
 say $\lambda_1,\ldots,\lambda_n$,  and assume that no two of the $\lambda_i$'s
 are conjugate to each other. Let
 $$\Gamma'\subset \Gamma_\mathbb{C}$$
 be the eigenspace associated to these eigenvalues. The last condition then implies that
 $$\Gamma_\mathbb{C}=\Gamma'\oplus\overline{\Gamma'},$$
 so that we have a complex torus $T$ defined as
 $$T=\Gamma_\mathbb{C}/(\Gamma'\oplus\Gamma).$$
 Note that $\phi_\mathbb{C}$ leaves stable $\Gamma'$
 and $\Gamma$, hence induces an endomorphism, which we will call $\phi_T$, of $T$.
 The endomorphism  of  $H_1(T,\mathbb{Z})\cong \Gamma$
 induced by $\phi_T$ is equal to $\phi$.

 The characteristic polynomial $f$ of $\phi$ has $\mathbb{Z}$-coefficients, and can be chosen to be any
 normalized degree 2n polynomial with integer coefficients,  subject to the condition
 (P).
 \begin{rema} \label{remark}Assume conversely that $T$ is a
$n$-dimensional complex torus admitting an endomorphism
$\phi_T$, whose induced morphism
$$\phi_{T*}:H_1(T,\mathbb{Z})\rightarrow H_1(T,\mathbb{Z})$$
identifies to $\phi$, via an isomorphism $\Gamma\cong H_1(T,\mathbb{Z})$.
Then $T$ is one of the tori  constructed as above, for a choice of $\lambda_i$'s.
\end{rema}We have now the following lemma:
 \begin{lemm}\label{lemma} If $n\geq2$ and the polynomial $f$
  is Galois, that is the Galois group of its splitting field
acts as the symmetric group on  the roots of $f$ (the eigenvalues of $\phi$),
 the torus $T$ is not an abelian variety.
 \end{lemm}
\begin{rema} This could be proved using the Albert classification of endomorphism
algebras of abelian varieties (see \cite{debarre}, \cite{birklange}). However, the proof
is so easy that it is better to do it directly.
\end{rema}

{\bf Proof.} Consider the rational N\'eron-Severi group
$$NS_\mathbb{Q}(T)\subset H^2(T,\mathbb{Q}),$$
generated over $ \mathbb{Q}$ by Chern classes of holomorphic line bundles on $T$.
This subgroup is contained in the subspace
$H^{1,1}(T)\subset H^2(T,\mathbb{C})$ of classes
representable by a form of type $(1,1)$, and is stable under the
action of $\phi^*_T$ on $H^2(T,\mathbb{Q})$.

Now the eigenvalues of $\phi_T^*$ on $H^{1,0}(T)\cong {\overline{\Gamma'}}^*$ are
the $\overline{\lambda_i}$, so that the eigenvalues of $\phi_T^*$ acting on
$$H^{1,1}(T)\cong
{\Gamma'}^*\otimes\overline{{\Gamma'}}^*$$
are exactly the $\lambda_i\overline{\lambda_j}$ for any $i,\,j$.

The subspace $NS_\mathbb{Q}(T)\otimes \mathbb{C}$ of $H^{1,1}(T)$
is stable under $\phi_T^*$, hence must be the eigenspace
associated to a set $S$ of  eigenvalues of the form $\lambda_i\overline{\lambda_j}$.
Since it is defined over $\mathbb{Q}$, the Galois group of the splitting field of $f$
has to leave stable the set $S$.
But by assumption, the Galois group acts on  the set
$\{\lambda_1,\ldots,\lambda_n,\overline{\lambda_1},\ldots,\overline{\lambda_n}\}$
as the symmetric group.
Hence it acts transitively on the products of two distinct elements of this set.
Since $n\geq2$, we can choose $i\not=j$, and then
$\lambda_i\lambda_j$ does not belong to the set $S$. Since $S$ is stable under the
Galois group, it follows that $S$ must be empty, so that $NS(T)=0$ and, a fortiori $T$ is not an abelian
variety.

\cqfd

\begin{rema} What we actually proved  is the equality $NS(T)=0$, under the same assumptions.
This is to be compared with \cite{shioda}, where Shioda proves that certain automorphism groups
acting on surfaces in projective space $ \mathbb{P}^3$ force the surface to have Picard group
to be as small as possible, namely generated by the class of the line bundle
$ \mathcal{O}(1)$.
\end{rema}

\section{\label{2} Construction of examples}
Let $(T,\phi)$ be a complex torus with endomorphism as in the previous section.

Our complex manifold $X$ will be constructed as follows : the diagonal of $T\times T$
and the graph of $\phi$ meet transversally at isolated points
$x_1=(0,0),\ldots,x_N$. Similarly the graph of $\phi_T$ meets $T\times 0$ transversally at finitely many points
$x_1=(0,0),\,y_1,\ldots,\,y_M$. We first blow-up
the points $x_1,\ldots,\,x_N,\,y_1,\ldots,\,y_M$  in $T\times T$. The proper transforms of the four
subtori
$$T\times 0,\,0\times T,\,T_{diag},\,T_{graph},$$
of $T\times T$, where the two last ones are respectively the diagonal and the graph of $\phi_T$,
are smooth and do not meet anymore. So we can blow-up their union, to get a
smooth compact manifold, which is K\"ahler since $T\times T$ is.
The following  shows that under the assumptions
of Lemma \ref{lemma}, $X$ satisfies the conclusion of Theorem \ref{mainintro}.

\begin{theo} \label{even}Assume $n\geq2$ and the characteristic polynomial of $\phi$ is Galois.
 Let $X'$ be a compact K\"ahler manifold. Assume that
there exists a ring isomorphism
$$\gamma: H^*(X',\mathbb{Z})\cong H^*(X,\mathbb{Z}).$$
Then $X'$ is not projective.
\end{theo}
{\bf Proof.} The isomorphism $\gamma$ shows us that $dim\,X'=dim X=2n$. Introduce the complex tori
$$Alb\,X,\,Alb\,X'.$$
They are complex tori of dimension
$$b_1(X)/2=2n,\, resp. \,\,b_1(X')/2=2n$$
defined respectively as
$$H^0(X,\Omega_X)^*/H_1(X,\mathbb{Z}),\,H^0(X',\Omega_{X'})^*/H_1(X',\mathbb{Z}).$$
We have the Albanese map, which is holomorphic :
$$alb_X:X\rightarrow Alb\,X,\,alb_{X'}:X'\rightarrow Alb\,X',$$
defined by integrating holomorphic forms along paths, well defined up to closed paths :
$$x\mapsto \int_{x_0}^x\in H^0(X,\Omega_X)^*,$$
$x_0\in X$ being a given base-point,
and similarly for $X'$.
This map induces by pull-back an isomorphism
$$alb_X^*:H^1(Alb\,X,\mathbb{Z})\rightarrow H^1(X,\mathbb{Z}),$$
and similarly for $X'$. Since
$H^*(Alb\,X,\mathbb{Z})\cong \bigwedge^*H^1(Alb\,X,\mathbb{Z})$, and $alb_X^*$
is compatible with the cup-product, we can identify
\begin{eqnarray}
\label{pullback}
alb_X^*:H^*(Alb\,X,\mathbb{Z})\rightarrow H^*(X,\mathbb{Z})
\end{eqnarray}
with the natural map induced by cup-product :
$$\bigwedge^*H^1(X,\mathbb{Z})\rightarrow H^*(X,\mathbb{Z}),$$ and similarly for $X'$.

In the case of $X$, the Albanese map is (up to translation) the natural blow-down map,
hence it is birational. This implies that in top degree $4n$, the map
$$alb_X^*:H^{4n}(Alb\,X,\mathbb{Z})\rightarrow H^{4n}(X,\mathbb{Z})$$
is an isomorphism.
As noted above, this map identifies to the cup-product map in top degree. Via the isomorphism $\gamma$,
we conclude now that the map
$$alb_{X'}^*:H^{4n}(Alb\,X',\mathbb{Z})\rightarrow H^{4n}(X',\mathbb{Z})$$
is also an isomorphism, which means that the map
$alb_{X'}^*$ is of degree $1$, that is birational.

Consider now the Gysin maps $alb_{X*}, \,alb_{X'*}$ on the cohomology of degree
$2$. They are
morphisms of Hodge structures, which can be defined here as the duals
with respect to Poincar\'e duality
of the pull-back maps $alb_{X}^*, \,alb_{X'}^*$ in degree $4n-2$, because the
considered manifolds have torsion free cohomology. Because the Albanese maps are of degree
$1$ in our case, the Gysin maps satisfy
\begin{eqnarray}\label{equation}alb_{X*}\circ alb_X^*=Id,\,alb_{X'*}\circ alb_{X'}^*=Id.
\end{eqnarray}
Since $alb_{X'}$ is birational, it induces an isomorphism
$$alb_{X'}^*:H^0(Alb\,X',\Omega^2_{Alb \,X'})\cong H^0(X',\Omega^2_{X'}).$$
This means that the morphism of Hodge structures
$$alb_{X'}^*:H^2(Alb\,X',\mathbb{Z})\rightarrow H^2(X',\mathbb{Z})$$
induces an isomorphism on $H^{2,0}$.
Equation  (\ref{equation}) shows then that the Hodge structure on
$Ker\,alb_{X'*}$ has no $(2,0)$-part, that is, it is made of cohomology classes of type
$(1,1)$ in the Hodge decomposition of $X'$.
\begin{rema} \label{rem10oct}
It is only to get this point that we needed  the
integral cohomology ring. The remainder of the proof will work as well with rational
cohomology, replacing in the following argument
(polarized) integral Hodge structures of weight $1$, (or equivalently
(projective) complex tori), with
(polarized) rational Hodge structures of weight $1$, (or equivalently
isogeny classes of
(projective) complex tori).
\end{rema}
We claim now that $Ker\,alb_{X'*}$ is the image under $\gamma^{-1}$ of the subgroup
$Ker\,alb_{X*}$. Indeed, Poincar\'e duality is given by cup-product followed
by the identification  given by the orientation
$$H^{4n}(\cdot,\mathbb{Z})\cong \mathbb{Z}$$ for all the considered manifolds.
Hence, since $\gamma$ is compatible with the cup-product, it is compatible up to sign with Poincar\'e duality.
Now we have seen that $\gamma$ identifies the images
of the pull-back maps (\ref{pullback}) for $X$ and $X'$.
Since $Ker\,alb_{X'*}$ in degree $2$ is the orthogonal complement with respect to Poincar\'e duality
of the image of the pull-back map in degree $4n-2$, and similarly for $X$, the claim follows.

In conclusion, we have a set of integral Hodge classes of degree $2$, provided by
the $\gamma^{-1}(\alpha),\,\alpha\in Ker\,alb_{X*}$. Recalling the
 construction of $X$ by a sequence of blow-ups, we see that this last group is freely generated
 by the classes $[\Delta_\cdot]$ of the exceptional divisors
 $$\Delta_{x_i},\Delta_{y_j},\,\Delta_{T\times 0},\,\Delta_{0\times T},\,\Delta_{diag},\,\Delta_{graph},$$
 over the corresponding centers of blow-up $x_i,\,y_j,\,T\times 0,\,0\times T,\, T_{diag},\,T_{graph}$.
 Let us denote $\delta_\cdot:=\gamma^{-1}([\Delta_\cdot])$.

 Each of these classes acts via cup-product on the cohomology of $X'$, and since we
 proved they are Hodge classes, their action is via morphisms of Hodge structures, of bidegree
 $(1,1)$. Hence the kernels of their restrictions to
 $H^1(X', \mathbb{Z})$ provide sub-Hodge structures of $H^1(X',\mathbb{Z})$.
 On the other hand, since $\gamma$ is compatible with the cup-product, these groups are
 the images under $\gamma^{-1}$ of the corresponding subgroups
$$Ker\,\cup [\Delta_\cdot]:H^1(X,\mathbb{Z})\rightarrow H^3(X,\mathbb{Z}).$$

 We look now at the maps $\cup [\Delta_\cdot]$. Consider the general
 situation of a blow-up
 $$\tau:\widetilde{Y}\rightarrow Y$$
 of a compact complex manifold $Y$ along a complex submanifold $Z\subset Y$
 of codimension $\geq2$.
 Let $j:E\hookrightarrow \widetilde{Y}$
 be the inclusion of the exceptional divisor, and let $\tau_E:E\rightarrow Z$ be the restriction of
 $\tau $ to $E$. Then the
 map
 $$\cup [E]\circ \tau^*:H^*(Y,\mathbb{Z})\rightarrow H^{*+2}(\widetilde{Y},\mathbb{Z})$$
 can be written as $$j_{E*}\circ \tau_E^*\circ j_Z^*,$$
 where $j_Z$ is the inclusion of $Z$ in $Y$.
 It is known that the composite  $j_{E*}\circ  \tau_E^*$ is injective
 on the cohomology of $Z$
 (cf \cite{voisin}, 7.3.3), hence we conclude
 that
 $$Ker\,\cup [E]\circ \tau^*=Ker\,j_Z^*.$$
 We apply this to $X$ and $T\times T$. (We ignore here the initial blow-up of points, since
 it does not affect $H^1$ or $H^3$, hence does not enter  in this computation.)
 For the complex subtorus $T\times 0$  we conclude that the kernel
 of the morphism
 $$\cup[\Delta_{T\times 0}]: H^1(X,\mathbb{Z})\cong H^1(T\times T,\mathbb{Z})
 \rightarrow H^3(X,\mathbb{Z})$$
 is equal to the kernel of the restriction map :
 $$H^1(T\times T,\mathbb{Z})\rightarrow H^1(T\times 0,\mathbb{Z}).$$
 Now the left hand side is isomorphic to
 $H^1(T,\mathbb{Z})\oplus H^1(T,\mathbb{Z})$ and the restriction map is the first projection
 $pr_1$.
 Similarly,  the kernel of the map
$ \cup[\Delta_{0\times T}]$ is equal to the kernel of the second projection $pr_2$.
Finally,
the kernel of $\cup[\Delta_{diag}]$ is equal to the kernel of
$pr_1+pr_2$, while the kernel of
$\cup[\Delta_{graph}]$ is equal to the kernel of
$pr_1+ \phi_T^*\circ pr_2$, since the restriction map from
$H^1(T\times T,\mathbb{Z})$ to $H^1(T_{graph},\mathbb{Z})=H^1(T,\mathbb{Z})$ identifies to
$pr_1+ \phi_T^*\circ pr_2$.

Let us now come back to $X'$. The above shows that $H^1(X',\mathbb{Z})$ contains four
integral
sub-Hodge structures $L_1,\,L_2,\,L_3,\,L_4$, images via $\gamma^{-1}$ of the $4$
subgroups
\begin{eqnarray}
\label{LI}
Ker\,pr_1,\,Ker\,pr_2,\,Ker\,pr_1+pr_2,Ker\, pr_1+\phi_T^*\circ pr_2.
\end{eqnarray}

 Consider the complex torus
 (dual to the Albanese torus) $$Pic^0(X')=H^1(X',\mathbb{C})/(H^{1,0}(X')\oplus H^1(X',\mathbb{Z})).$$
 Any integral sub-Hodge structure of $H^1(X',\mathbb{Z})$ provides
 a corresponding complex subtorus of $Pic^0(X')$ in an obvious way, namely, if
 $$L\subset H^1(X',\mathbb{Z})$$
 is a primitive sublattice such that
 $$L_{\mathbb C}=L^{1,0}\oplus L^{0,1},$$
 where $L^{1,0}=L_{\mathbb C}\cap H^{1,0}(X')$, and $L^{0,1}$ is its complex conjugate, then
 $$T_L= L_{\mathbb C}/(L^{1,0}\oplus L)$$
 is a complex subtorus of $Pic^0X'$.

 Hence we get four subtori $T_{L_i}$ of $Pic^0(X')$, which satisfy
 the following conditions, (because they are satisfied by the corresponding
 sublattices) :
 \begin{enumerate}
 \item $ T_{L_1}\oplus T_{L_2}= Pic^0(X')$.
 \item $T_{L_3}$ is isomorphic to $T_{L_1}$ and $T_{L_2}$ via the two projections
 induced by the previous isomorphism. In particular $T_{L_1}$ and $T_{L_2}$ are
 both isomorphic
 to some torus $T'$, and $Pic^0(X')\cong T'\oplus T'$.

 \item $T_{L_4}$ is isomorphic to $T'$ via the second  projection
 $p_2$ induced by the last isomorphism.
  Hence it is the (transpose of the) graph of an endomorphism $\phi'_{T'}$ of $T'$.
  \end{enumerate}

We finally observe that the action of $\phi'_{T'}$ on the homology $H_1(T',{\mathbb Z})$
is determined by the position of the four sublattices $L_i$. Using the description
(\ref{LI}) of the $\gamma^{-1}(L_i)$'s, we find that it identifies
via $\gamma$ to $\phi_T^*$,  that is to the dual of our initial endomorphism
$\phi$.

So we proved   that the variety $X'$ satisfies the property  that its Picard torus
is a product $T'\times T'$ where $T'$ is a complex torus admitting an automorphism
which acts on $H_1(T',{\mathbb Z})$ as the dual of $\phi$. We are now in position to apply
Lemma \ref{lemma}, combined with Remark \ref{remark}, which says that $T'$ cannot be
projective. On the other hand, if $X'$ were a projective variety, its Picard variety would be also
projective (cf \cite{voisin}, 7.2.2). So $X'$ is not projective.
\cqfd

 Since the $X$'s above have any possible even dimension $\geq4$ (because they are
 birational to products $T\times T$, with $dim\,T\geq2$), Theorem \ref{even} concludes
 the proof of Theorem \ref{mainintro} for even dimensions. In order to deal with odd dimensions, we
 prove with very similar arguments the following result :
 \begin{prop} Let $X$ be a variety constructed as above, and let $F$ be an elliptic curve.
 Let $X'$
be a K\"ahler compact manifold, such that there exists a ring isomorphism
 $$\gamma:H^*(X',\mathbb{Z})\rightarrow H^*(X\times F,\mathbb{Z}).$$
 Then $X'$ is not projective.
 \end{prop}

{\bf Proof.} Exactly as in the previous proof, we show first that the Albanese
map $alb_{X'}$ of $X'$ is birational,  that the kernel
of the map
$$alb_{X'*}:H^2(X',\mathbb{Z})\rightarrow H^2(Alb\,X',\mathbb{Z})$$
is made of Hodge classes, and is equal to
$$\gamma^{-1}(Ker\,(alb_{X\times F*}:
H^2(X\times F,\mathbb{Z})\rightarrow H^2(Alb\,(X\times F),\mathbb{Z}))).$$
The group $Ker\,(alb_{X\times F*}:
H^2(X\times F,\mathbb{Z})\rightarrow H^2(Alb\,(X\times F),\mathbb{Z}))$ is generated by the classes
of the exceptional divisors of the blowing down map
$$\tau:X\times F\rightarrow T\times T\times F.$$
The exceptional divisors are over centers which either are  of the form
$point\times F$ or are proper transforms of  subtori isomorphic to $T\times F$.
For any exceptional divisor $\Delta$ over $point\times F$, its class $[\Delta]$ induces
the morphism of Hodge structure
$$\cup[\Delta]:H^1(X\times F,\mathbb{Z})\rightarrow H^3(X\times F,\mathbb{Z})$$
which has for kernel (see the previous proof)
$$H^1(T\times T,\mathbb{Z})\subset H^1(T\times T\times F,\mathbb{Z}).$$
Indeed, this is also
$$Ker\,H^1(T\times T\times F,\mathbb{Z})\rightarrow H^1(point\times F,\mathbb{Z}).$$
Let $L\subset H^1(X',\mathbb{Z})$ be the image of this last kernel
via $\gamma^{-1}$. Then $L=Ker\,\cup\delta$, where $\delta=\gamma^{-1}([\Delta])$.
Since $\delta$ is a Hodge class, $L$ is a sub Hodge structure of $H^1(X',\mathbb{Z})$.
It follows that $Pic^0(X')$ contains a complex  subtorus
$T_L$ corresponding to the sub Hodge structure
$L$.

Now we consider the four other exceptional divisors $\Delta_i$ over the
proper transforms of
$$0\times T\times F,\,T\times0\times F,\,T_{diag}\times F,\,T_{graph}\times F.$$
There are, as in the previous proof, $4$ corresponding sub Hodge structures
$L_i$
of $H^1(X',\mathbb{Z})$ which are the kernels of
$\cup\delta_i,\,\delta_i=\gamma^{-1}([D_i])$.
It is immediate to see that the $L_i$ are contained in $L$, thus provide subtori
$T_{L_i}$ of $T_L$. As in the previous proof, one shows that :
\begin{enumerate}
 \item $ T_{L_1}\oplus T_{L_2}= T_L$.
 \item $T_{L_3}$ is isomorphic to $T_{L_1}$ and $T_{L_2}$ via the two projections
 induced by the previous isomorphism. In particular $T_{L_1}$ and $T_{L_2}$ are
 both isomorphic
 to some torus $T'$, and $T_L\cong T'\oplus T'$.

 \item $T_{L_4}$ is isomorphic to $T'$ via the second  projection
 $p_2$ induced by the last isomorphism.
  Hence it is the (transpose of the) graph of an endomorphism $\phi'_{T'}$ of $T'$.
  \end{enumerate}
Since the action of the endomorphism $\phi'_{T'}$ on its homology group $H_1$
is determined by the position of the lattices $L_i$, we see as in the previous proof that it has to
identify to the dual of $\phi$.
 In conclusion $T_L$ has to be a product $T'\times T'$, where $T'$ admits an endomorphism
which acts on its homology as the dual of $\phi$. Hence the subtorus
$T_L$ of $Pic^0(X')$ cannot be projective, so $Pic^0X'$ is not projective and $X'$ is not projective.
\cqfd
\section{Variants \label{3}}
\subsection{Other coefficients}
This subsection is due to P. Deligne (letter to the author).
\begin{theo} \label{thdelQ}
Let $X$ be as in Theorem \ref{even}. If $X'$ is
such that there exists a graded isomorphism
of rational cohomology rings
$$\gamma:H^*(X',\mathbb{Q})\cong H^*(X,\mathbb{Q}),$$
then $X'$ is not a  complex projective manifold.
\end{theo}
\begin{rema} Here, and in the previous section, we could conclude more generally that $X'$ cannot be smooth
complete algebraic or, in a more analytic language, Mo\"\i shezon.
 This is indeed sufficient to imply the existence of
   Hodge structures on the cohomology groups of $X'$, and
   of a polarized Hodge structure on $H^1(X')$,
which is  equivalent to the fact that $Pic^0(X')$ is a projective complex torus.
\end{rema}
The proof of Theorem \ref{thdelQ} is based on the following Lemma
\ref{ledel}: let $A^*=\oplus_kA^k$ be the rational cohomology ring of
a K\"ahler compact manifold (or a smooth complex
complete algebraic variety) and let
$A^*_\mathbb{C}:=A^*\otimes\mathbb{C}$.
Let $Z\subset A^k_\mathbb{C}$ be an algebraic subset which is defined by
homogeneous equations expressed only using the ring structure on $A^*$. The examples we shall consider
here and in the next subsection are :
\begin{enumerate}
\item $Z=\{\alpha\in A^k_\mathbb{C}/\alpha^l=0\,\,{\rm in }\,\,A^{kl}_\mathbb{C}\}$,
where $l$ is a given integer.
\item $Z=\{\alpha\in A^k_\mathbb{C}/\cup\alpha:A^{l}_\mathbb{C}\rightarrow A^{k+l}_\mathbb{C}\,\,
{\rm is\,\,not\,\,injective }\}$,
where $l$ is a given integer.
\end{enumerate}
\begin{lemm}\label{ledel} Let $Z$ be as above, and let $Z_1$ be an irreducible component of $Z$.
Assume the $\mathbb{C}$-vector space $<Z_1>$ generated by $Z_1$ is defined over $\mathbb{Q}$,
that is $<Z_1>=B^k_\mathbb{Q}\otimes\mathbb{C}$ for some $B_\mathbb{Q}^k\subset A_\mathbb{Q}^k$.
Then $B_\mathbb{Q}^k$  is a rational sub Hodge structure of $A^k_\mathbb{Q}$.
\end{lemm}

{\bf Proof.} It suffices to show that $B_\mathbb{C}^k=Z_1$ is stable under the Hodge decomposition
of $A^k_\mathbb{C}$. The Hodge decomposition can be seen as the character decomposition of the
action of $\mathbb{C}^*$ on $A^*_\mathbb{C}$ given by
$$z\cdot \alpha=z^p\overline{z}^q\alpha,\,\alpha\in A^{p,q}.$$
So it suffices to show that $<Z_1>$ is stable under this $\mathbb{C}^*$-action.
But $Z$ is stable under this action, by its definition and
because the action is compatible with the cup-product :
$$z\cdot(\alpha\cup\beta)=z\cdot\alpha\cup z\cdot\beta.$$
Being an irreducible component of $Z$, $Z_1$
 is also stable under this $\mathbb{C}^*$-action, and so is $<Z_1>$.
\cqfd
{\bf Proof of Theorem \ref{thdelQ}.} Let $P$ be
 the subspace
 of $H^2(X',\mathbb{Q})$ defined as the orthogonal complement with respect to Poincar\'e duality of
$\bigwedge^{4n-2}H^1(X',\mathbb{Q})\subset H^{4n-2}(X',\mathbb{Q})$.
As we mentioned in Remark \ref{rem10oct},
we used integral coefficients
in the  proof of Theorem \ref{even} only in order to conclude that $P$ consisted of Hodge classes.
Let now $P_0\subset P$ be defined as
$$P_0=\{\alpha\in P/\cup\alpha:H^1(X',\mathbb{Q})\rightarrow H^3(X',\mathbb{Q})\,\,{\rm is\,\,zero}
\}.$$
By Lemma \ref{ledel}, this is a sub Hodge structure of $P$. It is the image via $\gamma^{-1}
$ of the space generated by classes of exceptional divisors over points. Furthermore, we have for $\alpha\in
P/P_0$ an induced cup-product map
$$\cup\alpha:H^1(X',\mathbb{Q})\rightarrow H^3(X',\mathbb{Q}).$$
Introduce now the algebraic subset of $(P/P_0)\otimes\mathbb{C}$ :
$$Z=\{\alpha\in (P/P_0)\otimes\mathbb{C},\,\,\,\cup\alpha:H^1(X',\mathbb{Q})\rightarrow H^3(X',\mathbb{Q})\,\,
{\rm is\,\,not\,\,injective }\}.$$
This algebraic subset is the image via $\gamma^{-1}$ of
the union of four 1-dimensional $\mathbb{C}$-vector spaces, in fact defined over $\mathbb{Q}$,  generated
respectively by the classes of the four exceptional divisors over the blown-up complex subtori of $T\times T$.
We can thus apply Deligne's Lemma to each of these lines, thus concluding
that the four classes that we have denoted $\delta_{\cdot}$ in the  proof of Theorem \ref{even},
projected in $P/P_0$,
are Hodge classes for the Hodge structure on $P/P_0$.

Having this, we conclude   as in  the proof of Theorem \ref{even}.

\cqfd
In the above proof, we needed rational coefficients, in order to make sure that
the four lines above were defined over $\mathbb{Q}$. Deligne constructs now an example where the complex
cohomology algebra suffices to imply this.
The construction works as follows : start with the previous $X$. The four exceptional divisors
$\Delta_\cdot$ dominating the four complex subtori
$T\times 0,\,0\times T,\,T_{diag},\,T_{graph}$ are of the form
$\widetilde{T}_\cdot\times \mathbb{P}^{n-1}$, where $\widetilde{T}_\cdot$
is obtained by blowing-up finitely many (depending on which subtorus we consider)
points of $T$. Indeed,  the normal bundle $N$ of each of these subtori $T_\cdot$
in $T\times T$
is trivial, and if $\tau:\widetilde{T}_\cdot\rightarrow T_\cdot$ is the proper transform
of $T_\cdot$ under the initial blowing-up $\widetilde{T\times T}\rightarrow T\times T$ of points,
the normal bundle
of $\widetilde{T}_\cdot$ in $\widetilde{T\times T}$ is
isomorphic to
$\tau^*N(E)$, where $E$ is the exceptional divisor of $\tau$.

 Let us blow-up one
subvariety of the form $\widetilde{T}_{0\times T}\times \alpha_1$ in $\Delta_{0\times T}$, then
two subvarieties of the form $\widetilde{T}_{diag}\times \beta_1,\,\widetilde{T}_{diag}\times\beta_2$ in
$\Delta_{diag}$, and three subvarieties of the form $\widetilde{T}_{graph}\times\gamma_1,\,
\widetilde{T}_{graph}\times\gamma_2,\,
\widetilde{T}_{graph}\times\gamma_3$ in
$\Delta_{graph}$.
We  have constructed this way a smooth K\"ahler compact manifold $X_1$.
\begin{theo} \label{thdelC}
 If $X'_1$ is such that there exists a graded isomorphism of complex cohomology algebras
$$\gamma:H^*(X'_1,\mathbb{C})\cong H^*(X_1,\mathbb{C}),$$
then $X'_1$ is not a  complex projective manifold.
\end{theo}
{\bf Proof.} We use the same notations as before and consider the multiplication map
$$\cup\alpha:H^1(X'_1,\mathbb{Q})\rightarrow H^3(X'_1,\mathbb{Q}),$$
for $\alpha\in P/P_0$.
Let us introduce as in the previous proof the algebraic subset
$Z$ of $(P/P_0)\otimes \mathbb{C}$ consisting of those $\alpha$
for which $\cup\alpha$ is not injective. On one hand it is defined over $\mathbb{Q}$, and
on the other hand, using the isomorphism $\gamma$,  it consists
of the union of four complex spaces, which are of respective dimensions $1,\,2,\,3,\,4$.
(The reason is that, for $X_1$, the class of the exceptional divisor
 $\Delta_{T\times 0}$ generates the first of these spaces, while
 the second one is generated by the classes of the two exceptional divisors over
 $0\times T$ and so on.)
 It follows that each of these subspaces is in fact defined over $\mathbb{Q}$.
 Hence we can apply Lemma \ref{ledel}, to conclude that they are sub Hodge structures
 of $P/P_0$. Now for each of these spaces $V$, there is only one common kernel
 of the map $\cup\alpha,\,\,\alpha\in V$,
 which has to be a sub Hodge structure of $H^1(X'_1,\mathbb{Q})$.
 So we have produced again the four sub Hodge structures of $H^1(X'_1,\mathbb{Q})$, and
 the respective positions of the associated subspaces
 of $H^1(X'_1,\mathbb{C})$ are the same as for $X_1$ via the isomorphism $\gamma$.
We then concludes  that
$H^1(X'_1,\mathbb{Q})$ splits into the sum of two copies of a rational Hodge structure
which admits an automorphism which over $ \mathbb{C}$ is conjugate to $^t\phi$.
But then this automorphism is conjugate to $^t\phi$ also over $\mathbb{Q}$ and
we can conclude as in the proof of Theorem \ref{even}.

\cqfd
Note finally that if $Z$ is any reasonable simply connected compact topological
space,
Theorem \ref{thdelQ} (resp. \ref{thdelC}) remains true when
$X$ is replaced with $X\times Z$, resp. $X_1$ is replaced with $X_1\times Z$ :
the rational cohomology ring, resp. complex cohomology algebra of
$X\times Z$, resp. $X_1\times Z$, is not the one of a complex projective or Mo\"\i shezon variety.

Indeed, the $H^1$ is not modified by taking product with $Z$, and we recover
the space $P$ introduced in the above proofs as the annihilator of
$\bigwedge^{4n-2}H^1$ with respect to the cup-product.
The rest of the proof goes as before, leading to the conclusion
that  $H^1(Y,\mathbb{Q})$ cannot be endowed with a polarized Hodge structure, for any
$Y$ having the same rational (resp. complex) cohomology algebra as $X\times Z$, resp. $X_1\times Z$.

\subsection{A simply connected example}
In this subsection, we construct {\it simply connected} compact K\"ahler manifolds which do not have the
homotopy type (in fact the rational cohomology ring) of a projective complex manifold.

We start again with a torus $T$ as in section \ref{1},  endowed with an endomorphism
$\phi_T$, but assume now that its dimension $n$ is
$\geq 3$. Let us introduce the generalized Kummer variety
$$K=\widetilde{T/\pm1},$$
that is the desingularization of the quotient of $T$ by the $-1$ involution, obtained by blowing-up
the points of order $2$. This is a simply connected compact K\"ahler manifold. Its cohomology
$H^2(K,\mathbb{Q})$ is  the direct sum of
$\bigwedge^2 H^1(T,\mathbb{Q})$ and of the space generated by the exceptional divisors.

In $K\times K$, let us blow-up the diagonal, and then the proper transform of the
graph of $\phi_K$. ($\phi_T$ descends to a rational self-map of $K$, which we denote $\phi_K$.
Note that $\phi_K$ is not holomorphic, because one point of $T$ which is not of $2$-torsion can be sent
by $\phi_T$
to a $2$-torsion point of $T$, producing an indeterminacy point of
$\phi_K$. However, it is easy to see that the graph of $\phi_K$ is smooth, isomorphic
to the blow-up of $K$ at these points.)

 This will be our variety $X_2$. We shall  denote by
$\tau:X_2\rightarrow K\times K$ the blowing-down map.

\begin{theo}\label{thX2} If $X'_2$ is such that
there exists a graded isomorphism of rational cohomology rings
$$\gamma:H^*(X'_2,\mathbb{Q})\cong H^*(X_2,\mathbb{Q}),$$
then $X'_2$ is not a  complex projective manifold.
\end{theo}
The proof will proceed in several lemmas.
\begin{lemm} Consider the  subspaces
$$A^2_{i}:=\gamma^{-1}(\tau^*\circ pr_i^*(\bigwedge^2H^1(T,\mathbb{Q}))),\,i=1,\,2$$
of $H^2(X'_2,\mathbb{Q})$. Then $A^2_{i}$ are rational sub Hodge structures
of $H^2(X'_2,\mathbb{Q})$.
\end{lemm}
{\bf Proof.} Let $Z'$ be the algebraic subset of
$H^2(X'_2,\mathbb{C})$ defined as
$$Z'=\{\alpha\in H^2(X'_2,\mathbb{C})/\alpha^2=0\}.$$
This $Z'$ is the image under $\gamma^{-1}$ of the corresponding subset
$Z$ of
$H^2(X_2,\mathbb{C})$. Using the fact that
 $X_2$ has been deduced from $(T/\pm1)\times (T/\pm1)$ by a sequence of blow-ups with centers
of codimension $\geq3$, one sees easily that
$Z=Z_1\cup Z_2$, where
$$Z_i=\{\tau^*\circ pr_i^*(\alpha)/\alpha\in \bigwedge^2H^1(T,\mathbb{C}),\,\alpha^2=0 \,{\rm in}
\bigwedge^4H^1(T,\mathbb{C})\}.$$
But $\bigwedge^2H^1(T,\mathbb{C})$ is generated by those $\alpha$ such that $\alpha^2=0$.
It follows that $<Z_i>=\tau^*\circ pr_i^*(\bigwedge^2H^1(T,\mathbb{C}))$.

Applying Lemma \ref{ledel} to $Z'_i:=\gamma^{-1}(Z_i)$ gives the result.

\cqfd
Let now $A^*$ be the subalgebra of $H^*(X'_2,\mathbb{Q})$  generated by
$A_1^2\oplus A_2^2$ and let $P\subset H^2(X'_2,\mathbb{Q})$
be the orthogonal complement with respect to Poincar\'e duality of
$A^{4n-2}$. This space $P$ is the image under $\gamma^{-1}$
of the subspace $E$ of $H^2(X_2,\mathbb{Q})$ generated by the classes of exceptional divisors of $X_2$ over
$(T/\pm1)\times (T/\pm1)$.
This last space contains the two divisors classes $[\Delta_{diag}]$ and
$[\Delta_{graph}]$, and the divisors classes $[\Delta_{x_i}\times K],\, [K\times\Delta_{x_i}]$
where the $x_i$'s are the $2$-torsion points of $T$.
As in the previous section, consider the algebraic subset
$$Z'=\{\alpha\in P_{\mathbb{}}/\cup\alpha:A^2\rightarrow H^4(X'_2,\mathbb{C})\,\,{\rm
is\,\,not\,\,injective}\}.$$
$Z'$ is the image under $\gamma^{-1}$ of the corresponding
subset $Z$ of $E$. Now it is easy to check that
$Z$ is the union of four vector spaces defined over $\mathbb{Q}$, namely
$$<[\Delta_{diag}]>,\,<[\Delta_{graph}]>,\,<[\Delta_{x_i}\times K],\,i=1,\ldots,2^{2n}>,\,
<[K\times\Delta_{x_i}],\,i=1,\ldots,2^{2n}>.$$
Applying Lemma \ref{ledel}, we conclude that
the classes
$$\delta_{diag}:=\gamma^{-1}([\Delta_{diag}]),\,\delta_{graph}=\gamma^{-1}([\Delta_{graph}])$$
are Hodge classes in $H^2(X'_2,\mathbb{Q})$.
The kernels of the cup-product maps
$$\cup\delta_\cdot:A^2\rightarrow H^4(X'_2,\mathbb{Q})$$
are thus rational sub Hodge structures of
$A^2_1\oplus A^2_2$.
Examining via $\gamma$ their position  in $H^2(X_2,{\mathbb{Q}})$, we conclude that
$A_1^2$ and $A_2^2$ are isomorphic as rational Hodge structures and that
this rational Hodge structure carries an automorphism which acts as
$\phi_T^*={\bigwedge^2} {^t\phi}$ on $\bigwedge^2 H^1(T,\mathbb{Q})$.

Next, considering the proof of Lemma \ref{lemma}, we see that it actually shows the
following :
\begin{lemm} Assume the $\mathbb{Q}$-vector space $\bigwedge^2 H^1(T,\mathbb{Q})$
is endowed with a Hodge structure which is preserved by ${\bigwedge^2}{^t\phi}$. Then
either this Hodge structure is trivial, that is contains only Hodge classes,
 or it has no non zero Hodge classes.
\end{lemm}

We have now
\begin{lemm} If the Hodge structure on $A_i^2$ is trivial, then $X'_2$ cannot be
projective.
\end{lemm}
{\bf Proof.} Recall that $A_i^2=\gamma^{-1}(\tau^*\circ pr_i^*(\bigwedge^2 H^1(T,\mathbb{Q})))$.
It follows that it contains a $\mathbb{Q}$-vector subspace $V$
of dimension $\geq 2$ such that, for any $\alpha\in V,\,\alpha^2=0$ in $H^4(X'_2,{\mathbb{Q}})$.
(Take for $V$ a space of the form $\gamma^{-1}(\tau^*\circ pr_i^*(\alpha\wedge H^1(T,\mathbb{Q})))$
for a non zero $\alpha\in H^1(T,\mathbb{Q})$.)

This $V$ is then isotropic for any intersection form
on $H^2(X'_2,{\mathbb{Q}})$ of the form
$$q_c(\alpha)=c^{2n-2}\alpha^2,$$
where $c\in H^2(X'_2,\mathbb{Q})$. If $X'_2$ were projective, for an ample class $c$,
 this intersection form would have, by the Hodge index theorem (\cite{voisin}, 6.3.2),
  only one positive sign, when restricted to
   the space of Hodge classes of degree $2$. That would contradict the existence of this $2$-dimensional
 isotropic $V$ contained in $A_i^2$ if
 the Hodge structure on  $A_i^2$ were trivial.

\cqfd
\begin{rema} The proof even shows that under the same assumptions, $X'_2$ could not be K\"ahler.
\end{rema}
The two  last lemmas together show that the sub Hodge structure $A^2$ does not contain
any Hodge class if $X'_2$ is projective.
But then any degree $2$ Hodge class on $X'_2$ must be contained
in $P$. On the other hand, we have the following:
\begin{lemm}\label{last} For any class $c\in P$, the intersection form
$q_c(\alpha)=c^{2n-2}\alpha^2$  on $H^2(X'_2,\mathbb{Q})$ vanishes on $A^2$.
\end{lemm}
This is proved using the fact that $P$ is the image under $\gamma^{-1}$ of the space generated by the
(exceptional) divisor classes of $X_2$, and by proving the analogous result on $X_2$. It is
obviously essential here to assume that $n\geq 3$.
\cqfd
The proof of Theorem \ref{thX2} is now concluded by contradiction. If $X'_2$ were projective,
any degree $2$ Hodge class on $X'_2$ should be contained
in $P$, and $P$ should contain an ample class $c$. But then, again by the Hodge index theorem,
 the form $q_c$ would not vanish on
$A^2$, contradicting Lemma \ref{last}.

\cqfd

\end{document}